# Bounds for Solutions of Cubic Diophantine Equations

Howard Kleiman

**I. Introduction.** In this note, we use a simple argument to show that we are able to obtain solutions in $\mathbb{Z} \times \mathbb{Z}$ for some defining equations over $\mathbb{Q}(y)$ of degree 3.

**II. Preliminaries.** $x^3 - p(y)x + q(y)$ defines an absolutely irreducible polynomial over $\mathbb{Q}(y)$ where $p(y), q(y) \in \mathbb{Z}[y]$. Its discriminant i is $D = -4p^3 - 27q^2$. Using solution by radicals [2], its real root is of the form

$$\frac{\sqrt[3]{\frac{-27q + 3\sqrt{-3D}}{2}} + \sqrt[3]{\frac{-27q - 3\sqrt{-3D}}{2}}}{3}$$

Baker obtained the first effective bound for all solutions of $w^2 = f(y)$ in $\mathbb{Z} \times \mathbb{Z}$ [1] provided that $f(y)$ has at least three simple roots. Currently, better bounds exist. Call the best one, $B$.

**Theorem 1.** Let

$$x^3 - p(y)x + q(y) = 0 \quad (1)$$

(a) If $p(y)$ is not congruent to zero modulo 3, given current knowledge, we cannot obtain upper bounds for any solution of (1) in $\mathbb{Z} \times \mathbb{Z}$.

(b) If $p(y) \equiv 0 \pmod{3}$, assume that $D(y)$ has at least three simple roots. Then we can obtain all solutions $(x_0, y_0)$ of (1) in $\mathbb{Z} \times \mathbb{Z}$ satisfying the conditions that $w_0^2 = -3D(y_0)$ where $w_0$ is a rational integer.

**Proof.**

(a) Since $p(y)$ is not congruent to zero modulo 3, $D(y) \equiv 2(p(y))^3 \pmod{3}$ cannot be divisible by 3 for any rational integer $y = y_0$. But therefore $w_0^2 = -3D(y_0)$ cannot be divisible by 9, implying that $w_0$ cannot be a rational integer. Thus, we cannot use $B$ as an upper bound for $y_0$.

(b) In this case, any solution in $\mathbb{Z} \times \mathbb{Z}$ of $w^2 = -3D(y)$ has the property that

$-3D(y_0) = (-3)^2(C^2) \Rightarrow D(y_0) = -3C^2$ where $C \in \mathbb{Z}$. We thus have two possibilities: (a) $f(x, y_0)$ is irreducible and defines a metacyclic field, $M$, over $\mathbb{Q}$ such that by adjoining $i\sqrt{3}$ to it, $f(x, y_0)$ defines a cyclic field over $\mathbb{Q}(i\sqrt{3})$; (b) $f(x, y_0)$ is reducible and thus has a root that is a rational integer. In either case, using elementary algebra, we can use the factors of $q(y_0)$ in $\mathbb{Z}$ to determine if a rational integer root, $x_0$, of $f(x, y_0)$ exists. We need only continue testing until $|y_0| \leq B$.

*Comment.* It is likely that in most cases, $w_0$ is not rational even though $(x_0, y_0) \in \mathbb{Z} \; X \; \mathbb{Z}$. In particular, let $f(x, y_0) = (x - x_0)d(x, y_0)$. Here $d(x, y_0) = 0$ has a field discriminant of the form $-3r$ where $r$ is a power-free integer not divisible by 3.